\documentclass{article}
\usepackage{amsmath}
\usepackage{amssymb}
\usepackage{amsthm}
\newtheorem{theorem}{Theorem}
\newtheorem{lemma}[theorem]{Lemma}
\newtheorem{corollary}[theorem]{Corollary}

\newtheorem{proposition}[theorem]{Proposition}
\newtheorem{desired result}[theorem]{Desired result}

\newcommand{\begd}{\begin{displaystyle}}
\newcommand{\gl}{\lambda}

\newcommand{\gG}{\Gamma}
\newcommand{\ga}{\alpha}
\newcommand{\gb}{\beta}
\newcommand{\gd}{\delta}
\newcommand{\gD}{\Delta}
\newcommand{\gs}{\sigma}
\newcommand{\tr}{\textrm{tr}\,}
\newcommand{\ve}{\varepsilon}
\newcommand{\om}{\omega}
\newcommand{\mbq}{\mathbb{Q}}
\newcommand{\mbr}{\mathbb{R}}
\newcommand{\mbz}{\mathbb{Z}}
\newcommand{\mbc}{\mathbb{C}}

\newcommand{\mb}[1]{\mathbb{#1}}
\newcommand{\ol}[1]{\overline{#1}}
\newcommand{\mc}[1]{\mathcal{#1}}

\newcommand{\nequiv}{\equiv\hspace{-.13in}/\;}

\title{Trace formulas and class number sums}
\author{Nathan Jones}
\date{}

\begin{document}
\maketitle
{\def\thefootnote{}
\footnote{\today. \ {\it Mathematics Subject Classification (2000)}.
11R29, 11F32.}}
\begin{abstract}
\noindent
We specialize the Eichler-Selberg trace formula to obtain trace formulas for the prime-to-level Hecke action on cusp forms for certain congruence groups of arbitrary level.  As a consequence, we determine the asymptotic in the prime $p$ of the number of (weighted) $SL_2(\mbz)$-conjugation orbits of $2\times 2$ matrices of determinant $p$ whose reductions modulo $N$ lie in a given conjugacy class in $GL_2(\mbz/N\mbz)$.  This generalizes an 1885 result of Hurwitz.
\end{abstract}

\section{Introduction}
In \cite{hurwitz}, Hurwitz writes down formulas for sums of Hurwitz class numbers $H(-\gD)$ as $\gD$ runs through quadratic progressions to a prime modulus $N$.  He also mentions that these formulas may be generalized to the case where the modulus is not prime.  This paper generalizes Hurwitz's result to an arbitrary modulus $N$.  First, we describe all of this more precisely.

For any negative discriminant $\gD$, recall the Hurwitz class number
\[
H(-\gD) = \sum_{f(x,y) \in \mc{Q}^+_\mbz(\gD) \,//\, SL_2(\mbz)} \frac{2}{|SL_2(\mbz)_{f(x,y)}|}.
\]
Here we are denoting by
\[
\mc{Q}^+_\mbz(\gD) = \{ f(x,y) = \ga x^2 + \gb xy + \gamma y^2 : \; (\ga,\gb,\gamma) \in \mbz_{>0}\times \mbz^2, \; \gb^2-4\ga\gamma = \gD \}
\]
the set of positive definite integral binary quadratic forms of discriminant $\gD$, by $\mc{Q}^+_\mbz(\gD) \,//\, SL_2(\mbz)$ its orbit space with respect to the classical $SL_2(\mbz)$-action
\begin{equation*} \label{action}
f \cdot \begin{pmatrix} a & b \\ c & d \end{pmatrix} (x,y) = f(ax+by,cx+dy),
\end{equation*}
and by
\[
SL_2(\mbz)_{f(x,y)} = \{ A \in SL_2(\mbz) : \; f \cdot A = f \}
\]
the stabilizer in $SL_2(\mbz)$ of the form $f(x,y)$.  In addition, $H(0)$ is defined to be $-1/12$ and $H(m) = 0$ when $m < 0$.

Hurwitz shows, for example, that if $N$ is prime, $n > 1$ is coprime to $N$, and $a$ is any integer modulo $N$ with the property that $a^2-4n$ is a quadratic nonresidue modulo $N$, then
\[
(N+1) \sum_{t \equiv a \mod N} H(4n-t^2) = 2 \gs(n) + h_1^{(a)}\psi_1(n) + h_2^{(a)}\psi_2(n) + \dots + h_\mu^{(a)} \psi_\mu(n),
\]
where $\gs(n)$ is the sum of the divisors of $n$.  The $h_i^{(a)}$'s are coefficients which do not depend on $n$ and the $\psi_i(n)$'s are the Fourier coefficients of the $q$-expansions of certain weight $2$ cusp forms for the modular curve $X(N)$.  Thus, if we apply the Ramanujan bound $|\psi_i(p)| \leq 2p^{1/2}$ \cite{deligne}, we obtain
\begin{equation} \label{pasymptotic}
\sum_{t \equiv a \mod N} H(4n-t^2) = \frac{2}{N+1}\gs(n) + O_{N,\ve}(n^{1/2+\ve}).
\end{equation}

Let us re-interpret this asymptotic.  Note that, by pairing the positive definite form $f(x,y)$ with the negative definite form $-f(x,y)$ we have
\[
H(-\gD) = \sum_{f(x,y) \in \mc{Q}_\mbz(\gD) \,//\, SL_2(\mbz)} \frac{1}{|SL_2(\mbz)_{f(x,y)}|},
\]
where the sum is now taken over the orbit space of the set of \emph{all} integral binary quadratic forms of discriminant $\gD$.  Given integers $t$ and $n$, denote by
\[
\mc{T}(t,n) = \{ A \in M_{2 \times 2}(\mbz) : \; \tr A = t, \; \det A = n \}.
\]
If $t$ and $n$ satisfy $t^2-4n = \gD$, then there is a bijection
\begin{equation} \label{bijection}
\mc{Q}_\mbz(\gD) \longleftrightarrow \mc{T}(t,n)
\end{equation}
in which
\[
\ga x^2 + \gb xy + \gamma y^2 \leftrightarrow \begin{pmatrix} \frac{t+\gb}{2} & -\gamma \\ \ga & \frac{t-\gb}{2} \end{pmatrix}.
\]
This bijection is a map of $SL_2(\mbz)$-sets, where $SL_2(\mbz)$ operates by conjugation on $\mc{T}(t,n)$.  Thus we may re-write the Hurwitz class number as
\[
H \left( -(t^2-4n) \right) = \sum_{A \in \mc{T}(t,n) \,//\, SL_2(\mbz)} \frac{1}{|SL_2(\mbz)_A|}
\]
where $\mc{T}(t,n) \,//\, SL_2(\mbz)$ denotes the set of $SL_2(\mbz)$-conjugation orbits in $\mc{T}(t,n)$.  In this paper we prove
\begin{theorem} \label{asymptoticthm}
Let $N \geq 1 $ be any integer level, $n \geq 1$ a non-square integer coprime to $N$ and $\mc{A} \subset GL_2(\mbz/N\mbz)$ any $SL_2(\mbz/N\mbz)$-conjugation orbit with
\[
\det \mc{A} \equiv n \mod N.
\]
Then,
\begin{equation*} \label{asymptotic} 
\sum_{A \in \mc{T}_\mc{A}^e(n) \,//\, SL_2(\mbz)} \frac{1}{|SL_2(\mbz)_A|}  = \frac{2|\mc{A}|}{|SL_2(\mbz/N\mbz)|} \gs(n) + O_\ve(N^4n^{1/2+\ve}),
\end{equation*}
where
\[
\mc{T}_\mc{A}^e(n) := \{ A \in M_{2 \times 2}(\mbz) : \; A \mod N \in \mc{A}, \, \det A = n \text{ and } (\tr A)^2 < 4n \}.
\]
\end{theorem}
Note that this theorem specializes to \eqref{pasymptotic} in the case where $N$ is prime and $\mc{A}$ is the $SL_2(\mbz/N\mbz)$-conjugation orbit of trace $a$ and determinant $n$.

The case where $n=p$ is prime is of particular interst.  The work of Deuring \cite{deuring} (see also Theorem 14.18 of \cite{cox}) interprets the left-hand side of \eqref{pasymptotic} as essentially counting the number of isomorphism classes of elliptic curves over $\mbz/p\mbz$ whose frobenius endomorphism has trace congruent to $a$ modulo $N$.  Duke \cite{duke} uses this observation to unconditionally bound the mean-square error in the Chebotarev density theorem for the $N$-th division fields of elliptic curves over $\mbq$, for $N$ prime.  In a forthcoming paper we will use Theorem \ref{asymptoticthm} to strengthen Theorem $2$ of \cite{duke}.

\section{Acknowledgments}

This paper comprises a portion of my Ph. D. dissertation.  I would like to express gratitude to my advisor William Duke for his guidance.

\section{Statement of Results} \label{results}

Fix an arbitrary integer level $N \geq 1$.  In section \ref{groups} we define explicitly a family of subgroups
\[
\mc{D}_X \subset GL_2(\mbz/N\mbz),
\]
but for now we content ourselves with listing their relevant properties.  In order to obtain Theorem \ref{asymptoticthm} using trace formulas associated to the groups $\mc{D}_X$, we will make use of the following facts:
\begin{enumerate}
\item \label{capture}
The groups $\mc{D}_X$ capture all $SL_2(\mbz/N\mbz)$-conjugation orbits in $GL_2(\mbz/N\mbz)$, i.e.
\[
\forall A \in GL_2(\mbz/N\mbz), \; \exists \gamma \in SL_2(\mbz/N\mbz) \; \textrm{ and } \; \exists X \; \textrm{ so that } \gamma A \gamma^{-1} \in \mc{D}_X.
\]
\item \label{abelian}
Each group $\mc{D}_X$ is abelian, so that its space of class functions is spanned by its multiplicative characters $\chi$.
\item \label{technical}
For each $\mc{D}_X$ there exists a matrix $g_X \in GL_2(\mbr)$ so that $\det g_X = -1$ and, whenever $A \in M_{2 \times 2}(\mbz)$ satisfies $A \mod N \in \mc{D}_X$, we have
\[
g_X A g_X^{-1} \in M_{2 \times 2}(\mbz) \quad \text{ and } \quad g_X A g_X^{-1} \equiv A \mod N.
\]
\end{enumerate}
For now we leave the groups $\mc{D} = \mc{D}_X$ abstract, subject only to the above three conditions.  Our main formula is a trace formula for the action of $T_\mc{D}(n)$, the associated degree $n$ Hecke operator, on the space $S_k(\gG_\mc{D},\chi)$ of cusps forms with character $\chi$ relative to the associated congruence group $\gG = \gG_\mc{D}$ (for definitions, see Section \ref{background}).  In its statement we use the refined Hurwitz class number $H_\mc{C}(-\gD)$, which we define for any $SL_2(\mbz/N\mbz)$-conjugation orbit $\mc{C} \subseteq GL_2(\mbz/N\mbz)$ and any negative discriminant $t^2-4n$ by
\[
H_{\mc{C}}\left( -(t^2-4n) \right) := \sum_{{\begin{substack} {A \in \mc{T}(t,n) \,//\, SL_2(\mbz) \\ A \mod N \in \mc{C}} \end{substack}}} \frac{1}{|SL_2(\mbz)_A|}.
\]
Because of \eqref{bijection}, we see that $H_\mc{C}(-\gD)$ only depends on $\gD = t^2-4n$ and not individually on $t$ and $n$.  Whenever $t^2-4n > 0$ is a perfect square we define
\[
\mc{U}_\mc{C}(\sqrt{t^2-4n}) := \sum_{{\begin{substack} {A \in \mc{T}(t,n) \,//\, SL_2(\mbz) \\ A \mod N \in \mc{C}} \end{substack}}} 1.
\]
(One can show that $\sum_{\mc{C}} \mc{U}_\mc{C}(\sqrt{t^2-4n}) = \sqrt{t^2-4n}$.)  Also
\[
S(\mc{C},\chi) := \sum_{x \in \mc{C} \cap \mc{D}} \ol{\chi}(x)
\]
denotes the character sum and for $x \in M_{2 \times 2}(\mbz/N\mbz)$,
\[
\gG_x := \{ \gamma \in SL_2(\mbz) : \; ( \gamma \mod N ) x = x ( \gamma \mod N ) \}.
\]
We assume that
\begin{equation} \label{negativeIinG}
-I \in \gG, 
\end{equation}
which isn't particularly important but happens to be the case for our choice of groups $\mc{D}_X$.
\begin{theorem} \label{tfnonexplicit}
Suppose that $n$ is not a square and $\gcd(n,N)=1$.  Then
\[
\tr(T_\mc{D}(n)) = -t_h - t_e + \gd(\chi,k) \cdot \ol{\chi}(x_n) \cdot \gs(n),
\]
where
\[
t_h = \sum_{{\begin{substack} {\mc{C} = SL_2(\mbz/N\mbz) x SL_2(\mbz/N\mbz)^{-1}} \end{substack}}} [\gG_x : \gG] \sum_{{\begin{substack} {0 < d < \sqrt{n} \\ d \mid n} \end{substack}}} \mc{U}_\mc{C}(n/d-d) \cdot S(\mc{C},\chi) \cdot \frac{d^{k-1}}{n/d-d}
\]
and
\[
t_e = \sum_{{\begin{substack} {\mc{C} = SL_2(\mbz/N\mbz) x SL_2(\mbz/N\mbz)^{-1}} \end{substack}}} \frac{[\gG_x : \gG]}{2} \sum_{{\begin{substack} {t \in \mbz \\ t^2 < 4n} \end{substack}}} c_{t,k} \cdot S(\mc{C},\chi) \cdot H_{\mc{C}}\left( -(t^2-4n) \right).
\]
In the outermost sums, $\mc{C}$ is taken to range over all the $SL_2(\mbz/N\mbz)$-conjugation orbits in $GL_2(\mbz/N\mbz)$.
In the sum defining $t_e$,
\[
c_{t,k} := \frac{\rho_t^{k-1} - \ol{\rho_t}^{k-1}}{\rho_t - \ol{\rho_t}},
\]
where $\rho_t$ denotes any complex number satisfying $\rho_t^2-t\rho_t+n=0$,
\begin{equation} \label{delta} 
\gd(\chi,k) := \begin{cases}
               	1 & \textrm{ if } k = 2 \textrm{ and } \chi |_{\gG} \equiv 1 \\
		0 & \textrm{ otherwise},
              \end{cases}
\end{equation}
and $x_n \in \mc{D}$ is any element satisfying $\det x_n \equiv n \mod N$.
\end{theorem}

We remark that one need not assume $\mc{D}$ to be abelian, although it is convenient to simplify the proof.  All that is really necessary is that the multiplicative characters on $\mc{D}$ distinguish the $SL_2(\mbz/N\mbz)$ conjugation orbits in $\mc{D}$.  For example, if
\[
\mc{A} \cap \left\{ \begin{pmatrix} * & * \\ 0 & * \end{pmatrix} \mod N \right\} \neq \emptyset,
\]
then one can use the trace formula for $\gG_0(N)$ with character as developed in \cite{hijikata} to prove Theorem \ref{asymptoticthm}.  Otherwise we must use other congruence groups.  Chen \cite{chen} has used trace formulas for groups similar to ours (in the case of prime level and trivial character) to deduce the existence of isogenies between the jacobians of certain modular curves.

\section{Notation and Background} \label{background}

Throughout this paper we use the standard notation:
\[
\gG(N) := \{ \gamma \in SL_2(\mbz) : \; \gamma \equiv \begin{pmatrix} 
                                                       	1 & 0 \\
							0 & 1
                                                      \end{pmatrix}
\mod N \}.
\]
In particular, $\gG(1)$ denotes the full modular group $SL_2(\mbz)$.  For any subset $S \subseteq M_{2 \times 2}(\mbz/N\mbz)$ we put
\[
\mc{T}_S = \{ A \in M_{2 \times 2}(\mbz) : \; A \mod N \in S \}.
\]
Further we define, for any integers $t$ and $n$,
\[
\mc{T}_S(n) = \{ A \in \mc{T}_S : \; \det A = n \} \quad \text{ and } \quad \mc{T}_S(t,n) = \{ A \in \mc{T}_S(n) : \; \tr A = t \}.
\]
We abbreviate $\mc{T} := \mc{T}_{M_{2 \times 2}(\mbz/N\mbz)}$, so that our previous notation $\mc{T}(t,n)$ is consistent.

If $X$ is any set of matrices stable by left (resp. right) multiplication by a group $G$ of matrices, we use the usual notation
\[
 G \,\backslash\, X \quad \left( \text{ resp. } \; X \,/\, G \right)
\]
to denote the left (resp. right) coset space, whereas $X \,//\, G$ denotes the space of conjugation orbits, if $G$ acts on $X$ by conjugation.  Also, $Z(G)$ denotes the center of the group $G$, and $I$ denotes the $2\times 2$ identity matrix.

We now briefly set up the background, following \cite{miyake}, where full details may be found.  For a function $f$ on the upper half-plane we denote
\[
\left( f |_k \begin{pmatrix} a & b \\ c & d \end{pmatrix} \right) (z) := (ad-bc)^{k/2} (cz+d)^{-k} f \left(\frac{az+b}{cz+d}\right).
\]
Suppose $\gG$ is any Fuchsian group of the first kind and that
\[
\chi : \gG \longrightarrow \mbc^*
\]
is a multiplicative character whose kernel has finite index in $\gG$.  We consider the space of holomorphic weight $k$ modular forms with character $\chi$ for $\gG$
\[
\mc{M}_k(\gG, \chi) = \{ f : \mb{H} \rightarrow \mbc , f \textrm{ holomorphic (at cusps too), } \forall \gamma \in \gG, f|_k\gamma = \chi(\gamma)f \}.
\]
We note that whenever $-I \in \gG$ we have $\mc{M}_k(\gG, \chi) = \{0 \}$, unless 
\begin{equation} \label{chiassumption}
\chi(-I) = (-1)^k.
\end{equation}
Therefore we adopt this as a standing assumption.  The subspace of cusp forms is defined by
\[
\mc{S}_k(\gG, \chi) = \{ f \in \mc{M}_k(\gG, \chi) : f \equiv 0 \textrm{ at the cusps of } \gG \}.
\]
We recall the action of Hecke operators on these spaces.  Define the semigroup
\[
\tilde{\gG} := \{ g \in GL_2^+(\mbr) : \; [\gG : g \gG g^{-1} \cap \gG ] < \infty \; \textrm{ and } \;  [g \gG g^{-1}: g \gG g^{-1} \cap \gG ] < \infty \}.
\]
Let $\Upsilon$ be any subsemigroup satisfying
\[
\gG \subseteq \Upsilon \subseteq \tilde{\gG}
\]
and assume that $\chi$ extends to a multiplicative character of $\Upsilon$ so that for $\ga \in \Upsilon$ and $\gamma \in \gG$ we have
\begin{equation} \label{conditiononchi}
\ga \gamma \ga^{-1} \in \gG \Longrightarrow \chi(\ga \gamma \ga^{-1}) = \chi(\gamma).
\end{equation}
Then given any finite union of double cosets
\[
\mc{T} = \bigsqcup_{\ga \in \Upsilon} \gG \ga \gG
\]
we denote by $T$ (or by $T^\chi$, when we wish to emphasize the character $\chi$) the Hecke operator
\[
T : \mc{S}_k(\gG, \chi) \rightarrow \mc{S}_k(\gG, \chi),
\]
defined by the finite sum
\[
T(f) = \sum_{\ga \in \Upsilon} \det(\ga)^{k/2-1} \sum_{\ga_1 \in \gG \backslash \gG \ga \gG} \ol{\chi(\ga_1)} f|_k \ga_1.
\]
We refer to this situation by saying that the double coset space $\mc{T}$ \emph{defines} the Hecke operator $T$.

We now describe the specific Fuchsian groups and Hecke operators appearing in Theorem \ref{tfnonexplicit}.  Given the discussion in the previous paragraph, it remains to define $\gG$ and $\Upsilon$, explain which characters $\chi$ of $\gG$ we use and how they extend to $\Upsilon$, and finally to specify the double coset spaces $\mc{T}$ defining our Hecke operators.  Given any subgroup
\[
\mc{D} \subset GL_2(\mbz/N\mbz),
\]
we take
\[
\gG = \gG_\mc{D} := \mc{T}_\mc{D}(1) = \{ \gamma \in \gG(1) : \gamma \mod N \in \mc{D} \}
\]
and $\Upsilon$ to be the semigroup $\mc{T}_\mc{D}$.  We fix a group homomorphism
\begin{equation} \label{chionD}
\mc{D} \cap SL_2(\mbz/N\mbz) \longrightarrow \mbc^*.
\end{equation}
Precomposition with reduction modulo $N$ defines a character
\[
\chi : \gG_\mc{D} \twoheadrightarrow \mc{D} \cap SL_2(\mbz/N\mbz) \longrightarrow \mbc^*
\]
satisfying $\gG(N) \subseteq \ker \chi$.  Since $\mc{D}$ is assumed to be abelian, it is not difficult to show that any homomorphism \eqref{chionD} may be extended a homomorphism
\[
\mc{D} \rightarrow \mbc^*,
\]
In this way, $\chi$ extends to a semigroup homomorphism
\[
\chi : \mc{T}_\mc{D} \twoheadrightarrow \mc{D} \longrightarrow \mbc^*,
\]
and one verifies \eqref{conditiononchi} immediately.  We take as our Hecke operators $T = T_\mc{D}(n)$ to be those defined by the double coset space $\mc{T}_\mc{D}(n)$.

\section{A proof of Theorem \ref{tfnonexplicit}} \label{proofoftfnonexplicit}

We use the following more general trace formula due originally to Eichler \cite{eichler} (see also \cite{saito}, which works out the $\chi |_{\gG} =$ non-trivial case):
\begin{theorem} \label{traceformula}
Let $\gG$ be any Fuchsian group of the first kind, $\chi$ a character of $\gG$ of finite order and $k \geq 2$ an integer.  Assume that $\chi(-1) = (-1)^k$ if $-1 \in \gG$.  Let $T$ be any Hecke operator on the space $\mc{S}_k(\gG,\chi)$ of cusp forms for $\gG$ with character $\chi$.  Suppose that the double-coset space $\mc{T} \subset GL_2^+(\mbr)$ defining $T$ contains no scalar or parabolic elements.  Suppose further that there exists an element $g \in GL_2(\mbr)$ such that $\det(g)=-1$ and $\forall \ga \in \mc{T}$ one has
\begin{equation} \label{technicalintheorem}
g\ga g^{-1} \in \mc{T} \; \textrm{ and } \; \chi(g \ga g^{-1}) = \chi(\ga).
\end{equation}
Then the trace of $T$ is given by
\[
\tr(T) = - t_h - t_e + \gd(\chi,k) \sum_{\ga \in \gG \backslash \mc{T}} \ol{\chi}(\ga) ,
\]
where $\gd(\chi,k)$ is defined as in \eqref{delta},
\[
t_h = \frac{1}{|Z(\gG)|} \sum_{\ga \in \mc{T}^h \,// \, \gG} \ol{\chi(\ga)}\textrm{sgn}(\ga)^k\frac{\min\{|\eta_\ga|,|\zeta_\ga|\}^{k-1}}{|\eta_\ga -\zeta_\ga|}
\]
and
\[
t_e = \frac{1}{2} \sum_{\ga \in \mc{T}^e \,//\,\gG} \frac{\ol{\chi(\ga)}}{|\gG_\ga|}\frac{\eta_\ga^{k-1} - \zeta_\ga^{k-1}}{\eta_\ga - \zeta_\ga}.
\]
Here
\[
\mc{T}^h := \{ \ga \in \mc{T} : \; \tr(\ga)^2 > 4 \det(\ga) \textrm{ and $\ga$'s fixed points are cusps of } \gG \}
\]
and
\[
\mc{T}^e := \{ \ga \in \mc{T} : \; \tr(\ga)^2 < 4 \det(\ga) \}.
\]
The symbols $\eta_\ga$ and $\zeta_\ga$ refer to the complex eigenvalues of $\ga$, listed in either order.  When $\ga$ has real eigenvalues, $\textrm{sgn}(\ga)$ is defined to be the sign of either eigenvalue.
\end{theorem}
To deduce Theorem \ref{tfnonexplicit}, we apply Theorem \ref{traceformula} to the case $\gG = \gG_\mc{D}$ and $T = T_\mc{D}(n)$.  Since we assume $n$ is not a square, the double coset space $\mc{T} = \mc{T}_\mc{D}(n)$ doesn't have any scalar or parabolic elements.  The condition \eqref{technicalintheorem} follows immediately from the \ref{technical}rd condition describing the groups $\mc{D}_X$ in section \ref{results}.  Thus we may apply Theorem \ref{traceformula}.   For the hyperbolic term $t_h$ we fist observe that, using condition \eqref{chiassumption} to pair the summand belonging to $\ga$ with that belonging to $-\ga$, we have
\[
\sum_{\ga \in \mc{T}^h \,// \, \gG} \ol{\chi(\ga)}\textrm{sgn}(\ga)^k\frac{\min\{|\eta_\ga|,|\zeta_\ga|\}^{k-1}}{|\eta_\ga -\zeta_\ga|} = 2 \sum_{\ga \in \mc{T}^h_{>0} \,// \, \gG} \ol{\chi(\ga)} \frac{\min\{|\eta_\ga|,|\zeta_\ga|\}^{k-1}}{|\eta_\ga -\zeta_\ga|},
\]
where $\mc{T}^h_{>0} = \{ \ga \in \mc{T}^h : \; \textrm{sgn}(\ga) > 0 \}$.  Setting
\[
d = \min \{ |\eta_\ga|, |\zeta_\ga| \},
\]
we have $\{ |\eta_\ga|, |\zeta_\ga| \} = \{ d, n/d \}$.  We sort the matrices $\ga \in \mc{T}^h_{>0}$ occurring in $t_h$'s summation according to $d$.  Sorting the matrices $\ga$ occurring in the summations of both $t_h$ and $t_e$ according to
\[
\mc{C} = SL_2(\mbz/N\mbz) \left( \ga \mod N \right) SL_2(\mbz/N\mbz)^{-1},
\]
and noting \eqref{negativeIinG}, we see that
\begin{equation} \label{th}
t_h = \sum_{{\begin{substack} {\mc{C} \\ \mc{C} \cap \mc{D} \neq \emptyset} \end{substack}}} \sum_{{\begin{substack} {0 < d < \sqrt{n} \\ d \mid n} \end{substack}}} \frac{d^{k-1}}{n/d-d} \sum_{{\begin{substack} {\gb \in \mc{T}(n/d+d,n) \,//\, \gG(1) \\ \gb \mod N \in \mc{C}} \end{substack}}} \Phi_h(\gb)
\end{equation}
and
\begin{equation} \label{te}
t_e = \frac{1}{2} \sum_{{\begin{substack} {\mc{C} \\ \mc{C} \cap \mc{D} \neq \emptyset} \end{substack}}} \sum_{{\begin{substack} {t \in \mbz \\ t^2 - 4n < 0} \end{substack}}} c_{t,k} \sum_{{\begin{substack} {\gb \in \mc{T}(t,n) \,//\, \gG(1) \\ \gb \mod N \in \mc{C}} \end{substack}}} \Phi_e(\gb),
\end{equation}
where
\[
\Phi_h(\gb) = \sum_{\ga \in \left( \gG(1)\gb\gG(1)^{-1} \cap \mc{T}_\mc{D} \right) \,//\, \gG} \ol{\chi}(\ga)
\]
and
\[
\Phi_e(\gb) = \sum_{\ga \in \left( \gG(1)\gb\gG(1)^{-1} \cap \mc{T}_\mc{D} \right) \,//\, \gG} \frac{\ol{\chi}(\ga)}{|\gG_\ga|}.
\]
Now we evaluate $\Phi_h(\gb)$ and $\Phi_e(\gb)$.  Since $\mc{C} \cap \mc{D} \neq \emptyset$, we may as well assume that $\gb \mod N \in \mc{C} \cap \mc{D}$.
\begin{lemma}
Suppose that $\gb \equiv x \mod N \in \mc{D}$.  If $\gb$ is a hyperbolic matrix then
\begin{equation*} \label{hyperbolic}
\Phi_h(\gb) = [\gG_x : \gG] S(\mc{C},\chi),
\end{equation*}
while if $\gb$ is elliptic then
\begin{equation*} \label{elliptic}
\Phi_e(\gb) = \frac{[\gG_x : \gG]}{|\gG(1)_\gb|} S(\mc{C},\chi).
\end{equation*}
\end{lemma}
\begin{proof}
Consider reduction modulo $N$, which is a surjection:
\[
\textrm{red}_N : \gG(1)\gb\gG(1)^{-1} \cap \mc{T}_\mc{D} \twoheadrightarrow \mc{C} \cap \mc{D}
\]
Given $x \in \mc{C} \cap \mc{D}$, let $\ga_x \in \textrm{red}_N^{-1}(x)$ be arbitrary.  By definition we then have
\[
\textrm{red}_N^{-1}(x) = \gG_x \ga_x \gG_x^{-1}.
\]
Since $\mc{D}$ is abelian, the fiber $\textrm{red}_N^{-1}(x)$ is stable by $\gG$-conjugation.  Thus we may write
\[
\sum_{\ga \in \left( \gG(1)\gb\gG(1)^{-1} \cap \mc{T}_\mc{D} \right) \,//\, \gG} = \sum_{x \in \mc{C} \cap \mc{D}} \sum_{\ga \in \left( \gG_x \ga_x \gG_x^{-1} \right) \,//\, \gG}.
\]
Into how many $\gG$-conjugation orbits does $\gG_x \ga_x \gG_x^{-1}$ decompose?  Writing a right coset decomposition
\[
\gG_x = \bigsqcup_{b \in B} \gG b,
\]
we have
\begin{equation} \label{union}
\gG_x \ga_x \gG_x^{-1} = \bigcup_{b \in B} \gG b \ga_x b^{-1} \gG^{-1}.
\end{equation}
If $\gb$ is hyperbolic then so is $\ga_x$ and thus $\gG(1)_{\ga_x} = \{ \pm I \}$.  It follows that the above union is disjoint, and we obtain the lemma in this case.

For the elliptic case, if the centralizer $\gG(1)_{\gb} = \{\pm I \}$ then each $\gG(1)_{\ga_x} = \{ \pm I \}$, so that \eqref{union} is disjoint and there are exactly $[\gG_x : \gG ]$ $\gG$-conjugation orbits in $\gG_x \ga_x \gG_x^{-1}$.  In this case,
\[
\sum_{\ga \in \left(\gG_x \ga_x \gG_x^{-1} \right)//\,\gG} \frac{\ol{\chi}(\ga)}{|\gG_\ga|} = [\gG_x : \gG] \frac{\ol{\chi}(x)}{2}.
\]
Otherwise $\gG(1)_{\gb}$ is a group of order $4$ or $6$, and in that case we decompose the set $B$ of coset representatives into two subsets
\[
B = B_1 \sqcup B_2,
\]
where
\[
B_1 = \{ b \in B : \gG(1)_{b \ga_x b^{-1}} \subseteq \gG \}
\]
and
\[
B_2 = \{ b \in B : \gG(1)_{b \ga_x b^{-1}} \nsubseteq \gG \}
\]
and note that, for $b \in B_2$, $\gG(1)_{b \ga b^{-1}} \cap \gG = \{ \pm I \}$.  We then observe that, for any $b, b' \in \gG_x$ we have
\[
\gG b \ga_x b^{-1} \gG^{-1} = \gG b' \ga_x b'^{-1} \gG^{-1} 
\]
if and only if the equivalent conditions
\[
b'b^{-1} \in \gG(1)_{b'\ga_x(b')^{-1}} \gG \Longleftrightarrow b' \in \gG \gG(1)_{b \ga_x b^{-1}} b
\]
hold.  The first condition shows that unless $b, b' \in B_2$ we must have
\[
\gG b \ga_x b^{-1} \gG^{-1} \cap \gG b' \ga_x b'^{-1} \gG^{-1} = \emptyset,
\]
and when $b, b' \in B_2$ the second condition shows that the number of conjugation orbits in
\[
\bigcup_{b \in B_2} \gG b \ga_x b^{-1} \gG^{-1}
\]
collapses by a factor of $\frac{2}{|\gG(1)_{\gb}|}$.  In this case we have
\[
\begin{split}
\Phi_e(\gb) &= \sum_{x \in \mc{C} \cap \mc{D}} \left( \sum_{b \in B_1} \frac{\ol{\chi}(x)}{|\gG(1)_\gb|} + \frac{2}{|\gG(1)_\gb|} \sum_{b' \in B_2} \frac{\ol{\chi}(x)}{2} \right) \\
&= \frac{[\gG_x : \gG]}{|\gG(1)_{\gb}|} \sum_{x \in \mc{C} \cap \mc{D}} \ol{\chi}(x).
\end{split}
\]
\end{proof}
Returning to the proof of Theorem \ref{tfnonexplicit}, we see that, inserting our formulas for $\Phi_e(\gb)$ and $\Phi_h(\gb)$ into \eqref{te} and \eqref{th}, we obtain the expression for $t_e$ and $t_h$ as stated.

We finally show that for the remaining term we have
\[
\gd(\chi,k) \cdot \sum_{\ga \in \gG \backslash \mc{T}_\mc{D}(n)} \ol{\chi}(\ga) = \gd(\chi,k) \cdot \ol{\chi}(x_n) \cdot |\gG(1) \backslash \mc{T}(n)|.
\]
Both sides of the equation are nonzero if and only if $\chi$ restricted to $\gG$ is trivial.  In this case, all the $\chi(\ga)$'s are equal to any given $\chi(x_n)$.  We factor this out and then use the bijection of sets
\[
\gG_\mc{D} \,\backslash\, \mc{T}_\mc{D}(n) \longleftrightarrow \gG(1) \,\backslash\, \mc{T}(n), \quad \gG \ga \mapsto \gG(1) \ga.
\]
This finishes the proof of Theorem \ref{tfnonexplicit}.

\section{A proof of Theorem \ref{asymptoticthm}} \label{proofofasymptoticthm}

We proceed by taking appropriate linear combinations of the weight $k = 2$ formulas of Theorem \ref{tfnonexplicit}, making use of the character theory of abelian groups.  Because of \eqref{chiassumption}, we use the groups $\mc{D}/\{\pm I\}$ rather than $\mc{D}$.  In general let
\[
GL_2(\mbz/N\mbz) \longrightarrow GL_2(\mbz/N\mbz)/\{\pm I\} =: \ol{GL_2(\mbz/N\mbz)}, \quad \pm x \mapsto \ol{x}
\]
be the natural projection.  We note that each $\chi$ which satisfies \eqref{chiassumption} descends to a well-defined character on $\mc{D}$, which we continue to denote by $\chi$.  Our initial goal is to re-express our trace formulas in terms of $\ol{\mc{D}}$.  
For each $SL_2(\mbz/N\mbz)$-conjugation orbit $\mc{C} = SL_2(\mbz/N\mbz) x SL_2(\mbz/N\mbz)^{-1}$, the preimage of $\ol{\mc{C}}$ is $\pm \mc{C}$.  We make the definition
\[
S(\ol{\mc{C}},\chi) := \sum_{\ol{x} \in \ol{\mc{C}} \cap \ol{\mc{D}}} \chi(\ol{x}),
\]
and verify immediately that
\[
S(\ol{\mc{C}},\chi) = \begin{cases}
                       	S(\mc{C},\chi) & \textrm{ if } \; \mc{C} \neq -\mc{C} \\
			\frac{1}{2}S(\mc{C},\chi) & \textrm{ if } \; \mc{C} = -\mc{C}.
                      \end{cases}
\]
Looking at Theorem \ref{tfnonexplicit}, we see that in $t_e$, the term corresponding to $\mc{C}$ is equal to the term corresponding to $-\mc{C}$.  Thus, grouping these two terms together, we verify that when $k=2$ we have
\[
t_e = \sum_{{\begin{substack} {\ol{\mc{C}} \cap \ol{\mc{D}} \neq \emptyset} \end{substack}}} [\gG_x : \gG] S(\ol{\mc{C}},\chi) \sum_{{\begin{substack} { t \in \mbz \\ t^2-4n < 0} \end{substack}}} H_\mc{C}\left( -(t^2-4n) \right),
\]
where $\mc{C}$ has been chosen arbitrarily so that $\mc{C} \mapsto \ol{\mc{C}}$ and similarly $x$ is any element chosen above any $\ol{x} \in \ol{\mc{C}} \cap \ol{\mc{D}}$.  Similarly, we see that
\[
t_h = \sum_{{\begin{substack} {\ol{\mc{C}} \cap \ol{\mc{D}} \neq \emptyset} \end{substack}}} [\gG_x : \gG] \sum_{{\begin{substack} {0 < d < \sqrt{n} \\ d \mid n } \end{substack}}} \mc{U}_{\pm \mc{C}} (n/d-d) \cdot S(\ol{\mc{C}},\chi) \cdot \frac{d}{n/d-d}.
\]

Now suppose that
\[
\mc{A} = SL_2(\mbz/N\mbz) a SL_2(\mbz/N\mbz)^{-1} \subset GL_2(\mbz/N\mbz)
\]
is any conjugation orbit and $\det a \equiv n \mod N$.  Choosing the group $\mc{D}$ so that $\mc{D} \cap \mc{A} \neq \emptyset$, we can then assume that $a \in \mc{D}$.  Now we compute
\[
\frac{1}{|\ol{\mc{D}}^*|} \sum_{\chi \in \ol{\mc{D}}^*} \chi(\ol{a}) \tr(T^\chi_\mc{D}(n)).
\]
Using the orthogonality relations
\[
\frac{1}{|\ol{\mc{D}}^*|} \sum_{\chi \in \ol{\mc{D}}^*} \chi(\ol{a}) S(\ol{\mc{C}},\chi) = \begin{cases}
                                                                                            	1 & \text{ if } \ol{\mc{A}} = \ol{\mc{C}} \\
												0 & \text{ otherwise,}
                                                                                           \end{cases}
\]
we find that $\frac{1}{|\ol{\mc{D}}^*|} \sum_{\chi \in \ol{\mc{D}}^*} \chi(\ol{a}) \tr(T^\chi_\mc{D}(n))$ is equal to
\[
\frac{|\{ \chi \in \ol{\mc{D}}^* : \; \chi |_{\ol{\mc{D}} \cap SL_2(N)} = 1 \} |}{| \ol{\mc{D}} |} \cdot \gs(n) = \frac{2}{[\gG : \gG(N)]} \cdot \gs(n)
\]
minus
\[
[\gG_a : \gG] \left( \sum_{{\begin{substack} { t \in \mbz \\ t^2-4n < 0} \end{substack}}} H_\mc{A}\left( -(t^2-4n) \right) + \sum_{{\begin{substack} {0 < d < \sqrt{n} \\ d \mid n } \end{substack}}} \mc{U}_{\pm \mc{A}}(n/d-d) \frac{d}{n/d-d} \right).
\]
On the other hand using
\[
S_k(\gG(N)) = \bigoplus_{\chi \in \left( \gG/\gG(N) \right)^*} S_k(\gG, \chi) \quad \text{ and } \quad | \mc{D} / (\mc{D} \cap SL_2(\mbz/N\mbz)) | \leq N,
\]
together with the Ramanujan bound for the Hecke eigenvalues, we obtain Theorem \ref{asymptoticthm} with error term
\[
O_\ve \left( N \cdot \frac{\text{ the genus of } X(N)}{[\gG_a : \gG ]} n^{1/2+\ve} \right) = O_\ve (N^4 n^{1/2+\ve}).
\]
We have also used the bound
\[
\sum_{{\begin{substack} {0 < d < \sqrt{n} \\ d \mid n } \end{substack}}} \mc{U}_{\pm \mc{A}}(n/d-d) \frac{d}{n/d-d} \leq \sum_{{\begin{substack} {0 < d < \sqrt{n} \\ d \mid n } \end{substack}}} d \leq \sqrt{n} \sum_{d \mid n} 1 = O_\ve(n^{1/2+\ve}).
\]
For the genus of $X(N)$, see Theorem $4.2.11$ of \cite{miyake}.  Note that in case $n=p$ is prime we obtain the sharper error term $O(N^4p^{1/2})$, with an absolute constant.

\section{The groups $\mc{D}_X$} \label{groups}

It remains to define the groups $\mc{D}_X$ and verify the three conditions mentioned in section \ref{results}.
Let $q$ be an arbitrary integer and $\eta$ any integer co-prime to $N$.  The data
\[
X = (q , \eta , \gd) \in \mbz \times \{ \eta \in \mbz : \; \gcd(\eta, N) = 1 \} \times \{ 0, 1 \}
\]
determines the subgroup $\mc{D}(q,\eta,\gd) \subset GL_2(\mbz/N\mbz)$ by
\begin{equation} \label{definitionofD}
\mc{D}(q,\eta,\gd) = \left\{ \begin{pmatrix} x & \eta y \\ q \eta^* y & x + \gd L y \end{pmatrix} \mod N\right\} \subset GL_2(\mbz/N\mbz).
\end{equation}
Here $L$ is any fixed integer satisfying
\[
L \equiv 1 \mod 2^{n_2} \; \textrm{ and } \; \forall \text{ odd } p \mid N, \; L \equiv 0 \mod p^{n_p} \quad \left( N =: \prod_p p^{n_p} \right),
\]
$\eta^*$ is any integer satisfying
\[
\eta \cdot \eta^* \equiv 1 \mod N,
\]
and $x$ and $y$ run over all possible pairs of integers for which the given matrix modulo $N$ lies in $GL_2(\mbz/N\mbz)$.  The presence of $L$ is a technical nuisance arising from the fact that, when $N$ is even, the trace of a matrix can be either even or odd.

One calculates that each $\mc{D}(q,\eta,\gd)$ is abelian, so that property \ref{abelian} is satisfied.  We now verify properties \ref{capture} and \ref{technical}.

\subsection{Conjugation orbits of $SL_2(\mbz/N\mbz)$}

What information about a matrix $A \in M_{2 \times 2}(\mbz/N\mbz)$ specifies its $SL_2(\mbz/N\mbz)$-conjugation orbit?  Since the author has not seen a complete answer to this question in the literature, we give one here.  This will then be used in section \ref{capturesection} to show that the groups $D(q,\eta,\gd)$ satisfy property \ref{capture}.

Let $F$ be any matrix modulo $N$ and write
\begin{equation} \label{asube}
 F = \gl_F I + M_F A_F,
\end{equation}
where $M_F$ is the largest divisor of $N$ for which $F$ is scalar modulo $M_F$, $\gl_F \in \{ 0, 1, \dots, M_F - 1 \}$, and $A_F \in M_{2 \times 2}(\mbz/(N/M_F)\mbz)$ is the unique matrix making the above equation valid.  We shall see that the data of
\[
M_F, \, \gl_F, \; \text{ and } \, (\tr A_F , \det A_F ) \in (\mbz/(N/M_F)\mbz)^2
\]
uniquely determines the $GL_2(\mbz/N\mbz)$-conjugation orbit of $F$.  We proceed to give a complete list of the $SL_2(\mbz/N\mbz)$-conjugation invariants of any such $GL_2(\mbz/N\mbz)$-orbit.  To fix notation for our discussion, let 
\[
K = \prod_{p \textrm{ prime}} p^{m_p}
\]
be any integer level (soon to be $N/M_F$) and $A \in M_{2 \times 2}(\mbz/K\mbz)$ any matrix (soon to be $A_F$).  We denote the discriminant of $A$ by
\[
\gD(A) := \left( \tr A \right)^2 - 4 \det A \in \mbz/K\mbz.
\]
When $p$ is an odd prime dividing $K$ we define
\[
\chi_p : \{ A \in M_{2\times 2}(\mbz/K\mbz) : \; \gD(A) \equiv 0 \mod p, \, A \textrm{ nonscalar mod } p \} \rightarrow \{\pm 1 \}
\]
by
\[
\chi_p \left( \begin{pmatrix} a & b \\ c & d \end{pmatrix} \right) = \left( \frac{-b}{p} \right) \textrm{ or } \left( \frac{c}{p} \right),
\] 
whichever value is nonzero.  Here $\left( \frac{\cdot}{p} \right)$ denotes the Legendre symbol against $p$, which we will sometimes alternately denote by $\chi_p$, abusing notation and hoping not to cause confusion.  One checks that since $p$ divides $\gD(A)$, our expression for $\chi_p$ is well-defined.

If $m_2 \geq 2$ then for any $T \in \mbz/2^{m_2}\mbz$, we have
\[
T^2 \in \begin{cases}
         	\mbz/2^{m_2+1}\mbz & \; \text{ if $T$ is odd } \\
		\mbz/2^{m_2+2}\mbz & \; \text{ if $T$ is even.}
        \end{cases}
\]
Thus, whenever $4$ divides $K$ we have
\begin{equation} \label{whereisgD}
\gD(A) \in \{ x \in \mbz/2K\mbz : \; 2 \nmid x \} \sqcup \{ x \in \mbz/4K\mbz : \; 4 \mid x \}.
\end{equation}
If $4$ divides $K$ and $\gD(A)$ is even, then there is another invariant to define.  We take the characters $\chi_4$ and $\chi_8$ of conductors $4$ and $8$ respectively, defined to be $\equiv 0$ on the even integers and otherwise by
\[
\chi_4(1 \mod 4) = 1, \quad \chi_4(3 \mod 4) = -1
\]
and
\[
\chi_8(1 \mod 8) = \chi_8(7 \mod 8) = 1, \quad \chi_8(3 \mod 8) = \chi_8(5 \mod 8) = -1.
\]
We then define (again abusing notation)
\[
\chi_4 : \{ A \in M_{2\times 2}(\mbz/K\mbz) : \; \frac{\gD(A)}{4} \equiv 0,3 \mod 4 , \, A \textrm{ nonscalar mod }2\} \rightarrow \{\pm 1 \}
\]
by
\[
\chi_4 \left( \begin{pmatrix} a & b \\ c & d \end{pmatrix} \right) = \chi_4(-b) \textrm{ or } \chi_4(c),
\] 
depending on whether $-b$ or $c$ is odd.  If $m_2 \geq 3$ then we define
\[
\chi_8 : \{ A \in M_{2\times 2}(\mbz/K\mbz) : \; \frac{\gD(A)}{4} \equiv 0,2 \mod 8 , \, A \textrm{ nonscalar mod }2\} \rightarrow \{\pm 1 \}
\]
and
\[
\chi_4 \chi_8 : \{ A \in M_{2\times 2}(\mbz/K\mbz) : \, \frac{\gD(A)}{4} \equiv 0,6 \, \text{ mod}\, 8 , \, A \textrm{ nonscalar mod }2\} \rightarrow \{\pm 1 \}
\]
in the same way by
\[
\chi_8 \left( \begin{pmatrix} a & b \\ c & d \end{pmatrix} \right) = \chi_8(-b) \textrm{ or } \chi_8(c) 
\]
and
\[
\chi_4 \chi_8 \left( \begin{pmatrix} a & b \\ c & d \end{pmatrix} \right) = \chi_4(-b) \chi_8(-b) \textrm{ or } \chi_4(c) \chi_8(c),
\]
respectively.  One checks that these are all well-defined.

Now for any $T, D \in \mbz/K\mbz$ we set $\gD = T^2-4D$ and define formally, for each odd prime $p$ dividing $K$, the set of ``characters''
\[
\textbf{Char}_\gD(p^{m_p}) := \begin{cases}
                  	\{ \chi_p \} & \textrm{ if } p \mid \gD \\
			\emptyset & \textrm{ otherwise.}
                 \end{cases}
\]
If $2$ divides $K$ we define $\textbf{Char}_\gD(2^{m_2})$ by
\[
\textbf{Char}_\gD(2^1) := \emptyset, \quad \textbf{Char}_\gD(2^2) := \begin{cases}
                                               	\{ \chi_4 \} & \textrm{ if } \gD/4 \equiv 0, 3 \mod 4 \\
						\emptyset & \textrm{ otherwise,}
                                              \end{cases}
\]
and by
\[
\textbf{Char}_\gD(2^{m_2}) := \begin{cases}
			\{ \chi_4 \} & \textrm{ if } \gD/4 \equiv 3,4 \textrm{ or } 7 \mod 8 \\
			\{ \chi_8 \} & \textrm{ if } \gD/4 \equiv 2 \mod 8 \\
			\{ \chi_4\chi_8 \} & \textrm{ if } \gD/4 \equiv 6 \mod 8 \\
			\{ \chi_4, \chi_8, \chi_4\chi_8 \} & \textrm{ if } \gD/4 \equiv 0 \mod 8 \\
			\emptyset & \textrm{ otherwise}
                        \end{cases}
\]
if $m_2 \geq 3$.  (Here the notation ``$\gD/4 \equiv (\cdot) \mod \Box$'' is being used as shorthand for
\[
 4 \mid \gD \quad \text{ and } \quad \gD/4 \equiv (\cdot) \mod \Box.)
\]
Finally we set
\[
\textbf{Char}_\gD(K) :=  \bigcup_{m_p \geq 1} \textbf{Char}_\gD(p^{m_p}).
\]
For any matrix $A \in M_{2\times 2}(\mbz/K\mbz)$ of trace $T$ and determinant $D$ which is nonscalar modulo each prime $p$ dividing $K$, the above discussion shows that there is an evaluation map
\[
\textbf{Char}_\gD(K) \rightarrow \{ \pm 1 \}, \quad \chi \mapsto \chi(A).
\]
One calculates 
\[
\begin{pmatrix} 1 & 1 \\ 0 & 1 \end{pmatrix} \begin{pmatrix} a & b \\ c & d \end{pmatrix} \begin{pmatrix} 1 & -1 \\ 0 & 1 \end{pmatrix} = \begin{pmatrix} a+c & b+d-a-c \\ c & d-c \end{pmatrix}
\]
and
\[
\begin{pmatrix} 1 & 0 \\ 1 & 1 \end{pmatrix} \begin{pmatrix} a & b \\ c & d \end{pmatrix} \begin{pmatrix} 1 & 0 \\ -1 & 1 \end{pmatrix} = \begin{pmatrix} a-b & b \\ c+a-b-d & d+b \end{pmatrix},
\]
which shows, since the matrices $\begin{pmatrix} 1 & 1 \\ 0 & 1 \end{pmatrix}$ and $\begin{pmatrix} 1 & 0 \\ 1 & 1 \end{pmatrix}$ generate $SL_2(\mbz/K\mbz)$, that whenever $\gamma \in SL_2(\mbz/K\mbz)$ and $\chi \in \textbf{Char}_{\gD(A)}(K)$ we have
\[
\chi(A) = \chi(\gamma A \gamma^{-1}).
\]
The converse is also true:
\begin{lemma} \label{AconjugatetoB}
Let $K$ be a positive integer and let $A, B \in M_{2\times 2}(\mbz/K\mbz)$.  Assume that for each prime $p$ dividing $K$, neither $A$ nor $B$ is scalar modulo $p$.  Then, $A$ and $B$ are $GL_2(\mbz/K\mbz)$-conjugate to one another if and only if
\[
\tr A = \tr B \quad \text{ and } \quad \det A = \det B.
\]
When this is the case, $A$ and $B$ are $SL_2(\mbz/K\mbz)$-conjugate to one another if and only if we additionally have
\[
\chi(A) = \chi(B) \quad \forall \chi \in \textbf{Char}_{\gD}(K),
\]
where $\gD$ is the common discriminant of $A$ and $B$.
\end{lemma}
\begin{proof}
Since one direction was already observed above, suppose that $A$ and $B$ satisfy
all of the above hypotheses.  We will show that they are both $SL_2(\mbz/K\mbz)$-conjugate
to the matrix
\[
\begin{pmatrix} 0 & -\varpi^{-1}D \\ \varpi & T \end{pmatrix},
\]
where $\varpi \in (\mbz/K\mbz)^*$ is an element satisfying
\begin{equation} \label{varpiconsistency}
\chi(A) = \chi(B) = \chi(\varpi) \quad \forall \, \chi \in \textbf{Char}_\gD(K).
\end{equation}
To this end let $v = \begin{pmatrix} x \\ y \end{pmatrix} \in
(\mbz/K\mbz)^2$ be a variable vector and notice that the linear
transformation on $(\mbz/K\mbz)^2$ given by left multiplication by
$A$ has the form
\[
\left[ L_A \right]_{\{ v, \varpi^{-1}Av \}} = \begin{pmatrix} 0 & -\varpi^{-1}D \\ \varpi & T
\end{pmatrix}
\]
when written with respect to the ordered basis $\{v,\varpi^{-1}Av\}$ of
$(\mbz/K\mbz)^2$.  This verifies the claim, provided that we
can find a vector $v \mod K$ so that the change of basis matrix
belongs to $SL_2(\mbz/K\mbz)$, i.e. that we can solve
\begin{equation} \label{modpowerofp}
\phi(x,y) := cx^2+(d-a)xy-by^2 \equiv \varpi \mod p^{m_p} \quad \quad \left( A = \begin{pmatrix} a & b \\ c & d \end{pmatrix} \right)
\end{equation}
for each prime power $p^{m_p}$ exactly dividing $K$.  First assume $p$ is odd.  If $p$ divides both $b$ and $c$ then the polynomial $\phi(z,1) - \varpi$
has a simple zero mod $p$ which can be lifted to a zero mod $p^{m_p}$ by Hensel's lemma.  Otherwise, say $p \nmid c$.  Completing
the square we see that the condition
\[
\left( \frac{\gD y^2 + 4c\varpi}{p} \right) = 1 \text{ or } 0
\]
is necessary and sufficient for the solution of \eqref{modpowerofp}.  If $p$ does not divide $\gD$ then Lemma \ref{charactersum} below implies that we may choose $y$ so that this condition holds, while if $p$ divides $\gD$, then \eqref{varpiconsistency} asserts that the condition is true.

Now assume that $p=2$.  If $\gD$ is odd, define the polynomial
\[
g(z) := \begin{cases}
         \phi(z,1) - \varpi & \text{ if } b \text{ is odd} \\
	 \phi(1,z) - \varpi & \text{ otherwise.}
        \end{cases}
\]
Then $g$ has a simple zero mod $2$ which can be lifted to a zero mod $2^{m_2}$ by Hensel's lemma.

If $\gD$ is even, then since $A$ is non-scalar mod $2$, at least one of $b$ and $c$ is odd; let it be $c$.  Applying
the substitution 
\[
\left( x',y' \right) = \left( cx + \frac{d-a}{2}y, y \right),
\]
we arrive at the congruence
\[
(x')^2 \equiv \frac{\gD}{4} (y')^2 + c \varpi \mod 2^{m_2}.
\]
Using the fact that for $m_2 \geq 3$,
\begin{equation} \label{squaresmod2}
\left( (\mbz/2^{m_2}\mbz)^* \right)^2 = \{ x \mod 2^{m_2} : \; x \equiv 1 \mod 8 \},
\end{equation}
and the analogous fact for $m_2=2$, one checks that \eqref{varpiconsistency} guarantees that this transformed congruence may be solved.  
For example, if $\gD/4 \equiv 2 \mod 8$ we see that the congruence may be solved if $\varpi \equiv \pm c \mod 8$, i.e. if
\[
\chi_8(\varpi) = \chi_8(c).
\]
The other cases are similar.  This proves that $A$ and $B$ are $SL_2(\mbz/N\mbz)$-conjugate to one another.

To see that one only needs the trace and determinant condition to specify the $GL_2(\mbz/K\mbz)$-conjugation orbit,
note that
\[
\begin{pmatrix}
1 & 0 \\
0 & \varpi^{-1}
\end{pmatrix}
\begin{pmatrix}
0 & -\varpi^{-1} D \\
\varpi & T
\end{pmatrix}
\begin{pmatrix}
1 & 0 \\
0 & \varpi^{-1}
\end{pmatrix}^{-1}
=
\begin{pmatrix}
0 & - D \\
1 & T
\end{pmatrix}.
\]
\end{proof}
During the proof we made use of
\begin{lemma} \label{charactersum}
Let $p$ be an odd prime number and $\tau$, $\kappa \in (\mbz/p\mbz)^*$.  Then,
\[
\sum_{{\begin{substack} {j \neq 0 \mod p \\ \left(\frac{j}{p}\right) = \pm 1} \end{substack}}} \left( \frac{\tau+j\kappa}{p} \right) = -\left(\left(\frac{\tau}{p}\right) \pm \left(\frac{\kappa}{p}\right)\right)/2.
\]
\end{lemma}
\begin{proof} 
Let $S_{\pm}(\tau,\kappa)$ denote the sum on the left-hand side.
One sees immediately that the complete character sum
\begin{equation} \label{immediate}
S_\pm(\tau,\kappa) + S_\mp(\tau,\kappa) + \left( \frac{\tau}{p} \right) = 0.
\end{equation}
Suppose $\upsilon \in (\mbz/p\mbz)^*$ is a square modulo $p$.  We see that
\[
S_\pm(\tau,\kappa) = \left(\frac{\upsilon}{p}\right) S_\pm(\tau,\kappa) = \sum_{\left(\frac{j}{p}\right)=\pm 1} \left( \frac{\upsilon\tau+j\upsilon\kappa}{p} \right) = S_\pm(\upsilon\tau,\upsilon\kappa) = S_\pm(\upsilon\tau,\kappa).
\]
Thus,
\[
S_\pm(\upsilon\tau,\kappa) = S_\pm(\tau,\kappa) = S_\pm(\tau,\upsilon\kappa).
\]
We see that, for fixed $\kappa$, $S_\pm(\tau,\kappa)$ takes on one of two values, depending on the quadratic character of $\tau$.  Fix $\om$ any non-square modulo $p$.  We see that
\[
\frac{p-1}{2}\cdot S_\pm(\upsilon,\kappa) + \frac{p-1}{2}\cdot S_\pm(\om,\kappa) \pm \frac{p-1}{2}\left(\frac{\kappa}{p}\right) = \sum_{\tau \mod p} S_\pm(\tau,\kappa) = 0.
\]
On the other hand,
\[
-S_\mp(\tau,\kappa) = \sum_{\left(\frac{j\om}{p}\right) = \pm 1} \left(\frac{\om\tau + j\om\kappa}{p}\right) = S_\pm(\om\tau,\kappa).
\]
Thus we obtain
\[
S_\pm(\tau,\kappa) - S_\mp(\tau,\kappa) \pm \left(\frac{\kappa}{p}\right) = 0.
\]
This together with \eqref{immediate} proves the lemma.
\end{proof}

Returning to our discussion of $SL_2(\mbz/N\mbz)$-orbits, we define formally, for $F$ as in \eqref{asube},
\[
\chi(F) := \chi(A_F), \quad \chi \in \textbf{Char}_{((\tr A_F)^2 - 4 \det A_F)}(N/M_F).
\]
Then putting $K = N/M_F$ in the previous lemma shows the following.
\begin{proposition} \label{EconjugatetoF}
A pair of matrices $F$ and $F'$ in $M_2(\mbz/N\mbz)$ are $GL_2(\mbz/N\mbz)$-conjugate to one another if and only if $M_F = M_{F'} =: M$, $\gl_F = \gl_{F'}$,
\[
\tr A_F \equiv \tr A_{F'} \mod N/M, \quad \text{ and } \quad \det A_F \equiv \det A_{F'} \mod N/M.
\]
When this is the case, $F$ and $F'$ are $SL_2(\mbz/N\mbz)$-conjugate to one another if and only if we additionally have
\[
\chi(F) \equiv \chi(F') \quad \forall \; \chi \in \textbf{Char}_{\gD}(N/M),
\]
where $\gD$ is the common discriminant $(\tr A_F)^2 - 4\det A_F = (\tr A_{F'})^2 - 4\det A_{F'}$.
\end{proposition}

The proposition shows that for any conjugacy class $\mc{C} \subset GL_2(\mbz/N\mbz)$, we have
\[
|\, \mc{C} \,//\, SL_2(\mbz/N\mbz)| \leq 2^{|\{ p : p \mid N \}|+2} \ll N.
\]
This implies
\begin{corollary}
Let $\mc{C} \subset GL_2(\mbz/N\mbz)$ be any conjugacy class and $p$ a prime with $p \equiv \det \mc{C} \mod N$.  Then
\[
\sum_{\ga \in \mc{T}_\mc{C}^e(p) \,//\, \gG(1)} \frac{1}{|\gG(1)_\ga|} = \frac{2 |\mc{C}|}{|SL_2(\mbz/N\mbz)|} p + O(N^5 p^{1/2}),
\]
with an absolute constant.
\end{corollary}
We also remark that one may write down a more explicit version of Theorem \ref{tfnonexplicit} by incorporating the information from Proposition \ref{EconjugatetoF}.  However, we will not do so here.

\subsection{The groups $\mc{D}_X$ satisfy property \ref{capture}} \label{capturesection}

We are now in a position to show
\begin{proposition} \label{captureprop}
The groups $D(q,\eta,\gd)$ capture the $SL_2(\mbz/N\mbz)$-conjugation orbits in $GL_2(\mbz/N\mbz)$.  
\end{proposition}
\begin{proof}
Let $F$ be any matrix in $GL_2(\mbz/N\mbz)$, written as in \eqref{asube}.  Replacing $x$ by $\gl_F + M_F x$ and $y$ by $M_F y$ in \eqref{definitionofD}, we have
\[
 \gl_F I + M_F \begin{pmatrix} x & \eta y \\ q\eta^* y & x+\gd Ly \end{pmatrix} \in D(q,\eta,\gd).
\]
Let $T = \tr A_F$ and $D = \det A_F$ and define the exponents $m_p$ by
\[
N/M_F =: \prod_p p^{m_p}.
\]
By Proposition \ref{EconjugatetoF}, such a matrix is $SL_2(\mbz/N\mbz)$-conjugate to $F$ if and only if for each prime $p$ dividing $N/M_F$, 
\begin{equation} \label{firstcondition}
T \equiv 2x + \gd L y \mod p^{m_p}, \quad D \equiv x^2 + \gd L xy - q y^2 \mod p^{m_p}, \quad p \nmid y 
\end{equation}
is solvable in $x$ and $y$ simultaneously with
\[
\forall \chi \in \textbf{Char}_\gD(p^{m_p}), \; \chi(\eta y) = \chi(F).
\]
By the Chinese Remainder Theorem, we can always choose $\eta$ so that this second condition is satisfied, so it suffices to solve  \eqref{firstcondition}.  For an odd prime $p$, the equations read
\[
T \equiv 2x \mod p^{m_p}, \quad D \equiv x^2 - q y^2 \mod p^{m_p}, \quad p \nmid y.
\]
If $T^2-4D \equiv 0 \mod p^{m_p}$ then these may be solved if and only if $p^{m_p}$ divides $q$.  Otherwise the equations may be solved if and only if
\[
p^{m_p'} \mid \mid q \quad \text{ and } \quad \frac{q}{p^{m_p'}} \mod p^{m_p-m_p'} \in \frac{T^2-4D}{p^{m_p'}} \left(\left( \mbz/p^{m_p-m_p'}\mbz \right)^*\right)^2,
\]
where $m_p' < m_p$ is the exponent of $p$ dividing $T^2-4D$.  This can be arranged by choosing $q$ appropriately modulo $p^{m_p}$.

If $p = 2$ then we put $\gd = 0$ if $T$ is even and $\gd = 1$ otherwise.  In this case \eqref{firstcondition} reads
\[
T \equiv 2x +\gd y \mod 2^{m_2}, \quad D \equiv x^2 + \gd xy - q y^2 \mod 2^{m_2}, \quad 2 \nmid y.
\]
If $m_2 = 1$ then this can easily be solved.  In case $m_2 > 1$, we note that if $T$ is even then $(T^2-4D)/4$ is well-defined modulo $2^{m_2}$ and so as in the case of odd $p$, we may find conditions on $q$ guaranteeing that \eqref{firstcondition} can be solved.  If $T$ is odd then the equations may be solved if and only if
\begin{equation} \label{1plusfourq}
1+4q \in (T^2-4D)\left( (\mbz/2^{m_2}\mbz)^* \right)^2.
\end{equation}
For the ``if'' direction, choose any $y$ modulo $2^{m_2}$ so that $(1+4q)y^2 \equiv T^2-4D \mod 2^{m_2}$ and choose any $x$ so that $T \equiv 2x+y$.  By \eqref{whereisgD} we have
\[
D \equiv x^2+xy-qy^2 \mod 2^{m_2-1}.
\]
In case $D \nequiv x^2+xy-qy^2 \mod 2^{m_2}$, replace $x$ by $x+2^{m_2-1}$.

Using \eqref{squaresmod2}, one sees that choosing $q \equiv D \mod 2$ guarantees \eqref{1plusfourq}.  One checks that this choice works for $m_2 = 2$ just as well.  This proves Proposition \ref{captureprop}.
\end{proof}

\subsection{The groups $\mc{D}_X$ satisfy property \ref{technical}}

\begin{proposition}
For each $\mc{D}(q,\eta,\gd)$ there exists a matrix $g \in GL_2(\mbr)$ so that $\det g = -1$ and, whenever $A \in M_{2 \times 2}(\mbz)$ satisfies $A \mod N \in \mc{D}(q,\eta,\gd)$, we have $g A g^{-1} \in M_{2 \times 2}(\mbz)$ and 
\[
g A g^{-1} \equiv A \mod N.
\]
\end{proposition}
\begin{proof}
Let $\gamma \in \gG(1)$ be any matrix satisfying
\[
\gamma \equiv \begin{pmatrix} \eta & 0 \\ \gd L & \eta^* \end{pmatrix} \mod N
\]
and
\[
\gamma_q = \frac{1}{\sqrt{q}} \begin{pmatrix} 0 & 1 \\ q & 0 \end{pmatrix}.
\]
We can take $g = \gamma \gamma_q$.
\end{proof}

\vspace{.1in}

\begin{center}
Centre de Recherches Math\'{e}matiques \\
Universit\'{e} de Montr\'{e}al \\
P.O. Box 6128, \\
Centre-ville Station \\
Montr\'{e}al, Qu\'{e}bec  H3C 3J7, Canada. \\
E-mail:  jones@dms.umontreal.ca
\end{center}
\end{document}